\documentclass[submission,copyright,creativecommons]{eptcs}
\providecommand{\event}{Post Proceedings of ADG'23} 
\usepackage[T1]{fontenc}
\usepackage{xurl,amsmath}
\usepackage{cite}
\usepackage{graphicx}
\usepackage{framed}
\usepackage{breqn}
\usepackage{amsfonts}
\usepackage{listings}
\usepackage{color}
\usepackage{subfigure}
\usepackage{multicol}
\usepackage{longtable}
\usepackage{booktabs}
\usepackage{amssymb}
\usepackage{doi}

\begin{document}
\title{Solving Some Geometry Problems of the Náboj 2023 Contest with Automated Deduction in GeoGebra Discovery}
\def\titlerunning{Solving Some Geometry Problems of the Náboj 2023 Contest with GeoGebra Discovery}
\def\authorrunning{A.~Hota et al.}
\author{
Amela Hota
\institute{The Private University College of Education of the Diocese of Linz, Austria}
\email{amela.hota@ph-linz.at}
\and
Zolt\'an Kov\'acs 
\institute{The Private University College of Education of the Diocese of Linz, Austria}
\email{zoltan.kovacs@ph-linz.at}
\and Alexander Vujic
\institute{The Private University College of Education of the Diocese of Linz, Austria}
\email{alexander.vujic@ph-linz.at}}
\maketitle              
\begin{abstract}
In this article, we solve some of the geometry problems of the Náboj
2023 competition with the help of a computer, using examples that
the software tool GeoGebra Discovery can calculate. In each case, the calculation requires
symbolic computations. We analyze the difficulty of feeding the problem
into the machine and set further goals to make the problems of this type
of contests even more tractable in the future.
\end{abstract}

\section{Introduction}

With the everyday rise of Artificial Intelligence (AI), the power of computers
has become tangible for the masses. Yes, it can do your homework (not
just in maths), but it can also pass your A-level exams.\footnote{%
See \url{https://telex.hu/tech/2023/05/09/chatgpt-bing-mesterseges-intelligencia-matematika-erettsegi}.}
A long series
of ad-hoc studies have shed light on what the present can offer: often
instant and perfect answers to questions that take years of learning to
solve by human means. This raises a number of research questions, such as whether
the current school system is still needed, whether teachers are still
needed, or whether it is enough to have AI.\footnote{%
See Bill Gates' notes at \url{https://www.gatesnotes.com/ASU-and-GSV?WT.mc_id=20230419100000_ASU-GSV-2023_BG-EM_}.}
Of course, alongside the praise, there are also many criticisms:
AI sometimes makes mistakes, especially in textual
problem settings where the question is formulated in a challenging way.

Automatic geometrical derivations, on the other hand, are perfect and, as
such, there is no such a major possibility of error. The answer is not derived from a
statistically computed result (as is so often the case with AI-based
algorithms), but a verifiable derivation is given in each case. We do not
claim that the two directions cannot meet once, and indeed, ultimately,
AI should refer to, i.e.~use, the ADG algorithm as a subroutine. There
are already prototypes working in this direction, e.g.~the WolframAlpha
system has been successfully coupled with an AI frontend.\footnote{%
See Stephen Wolfram's notes at \url{https://writings.stephenwolfram.com/2023/03/chatgpt-gets-its-wolfram-superpowers/}.}

In this contribution, we aim for less. We are just trying to solve competitive
problems with an ADG algorithm in the background. However, we leave the
exact task setting to the user. This means that it is up to the user to
provide the exact flow of the editing task with concrete steps. This must
be done in GeoGebra Discovery\footnote{GeoGebra Discovery is freely available at
\url{https://kovzol.github.io/geogebra-discovery}.}. However, for the inference, which requires
a symbolic calculation in the tasks, the ADG algorithm steps in and, as
we will see, gives the correct result in all cases. In the second half of
the paper, we propose how the range of problems that can be solved in
this way can be further extended.

\section{The Náboj Contest}

According to \href{www.naboj.org}{\textbf{naboj.org}}, Náboj is an international mathematical competition designed for teams of five high-school students that represent their schools, which lasts 120 minutes and where they are trying to solve as many given problems as possible. As soon as the team correctly solves any of the problems, they receive new ones. The solutions of the problems are usually numerical. The team that solves most problems correctly in the given time limit wins. The Náboj problems in contrast to the most school exercises require a certain amount of inventiveness and ingenuity.

Traditionally, many geometric tasks require proving. Checking if a proof is correct may be a difficult process for the organizers, so it is usually avoided to set proof related problems during contests like Náboj. Instead, problem settings require computing fractions, or better, providing a non-trivial algebraic number. As a consequence, geometric problems in Náboj are mostly non-geometric, or if still so, they are set in a way to require a numerical result.

All the problems we discuss in this paper will have exact answers, non-trivial fractions or some root expressions. It is clear that the exact
definition of the latter requires symbolic computation. By default,
software that allows a geometric problem to be well visualized (GeoGebra
in particular) provides only numerical support for measuring the quantity
in question. However, the software presented here, the GeoGebra fork
GeoGebra Discovery, is
capable of making measurements symbolically. This also means that a full
proof has already been created in the background, but the user is not
informed about this.

The use of electronic aids in the Náboj competition has recently been
restricted since the competition is on-site again. In the long
term, electronic assistance will certainly not be eliminated. It is a fact that students are turning to AI for quick help, and the tasks set must take this into account. It may seem like fun, but the rapid pace of the modern age
also poses a huge challenge for assignment writers: is the task set
difficult enough to prevent the AI from giving a quick, accurate,
complete answer? But the task setter is not only fighting the AI, but
also the ADG algorithm: cannot the task be solved in a flash if the right
data is entered into the right software in the right way and the right
button is pressed?

Overall, we conclude that the tasks set should be tested with different
software before they are announced, in order to avoid embarrassing
surprises. Even if we manage to keep the students working with paper and
pencil only for the 2 hours of the competition, i.e.~to exclude
electronic assistance, it raises serious questions about what the AI and
ADG algorithms can achieve with the tasks set. In the preliminary
analysis of problems, but also in retrospect, when we return to the
correct solution of problems in mathematics class or in a specialised
course, it may be useful to use the electronic method.

\section{Mathematical Background}

The method used by GeoGebra Discovery is essentially the Recio-Vélez
method \cite{RecioVelez99}, complemented by the algorithm given in \cite{scsc-2020}. For the problems of the
Náboj competition we are discovering a ratio of two lengths, it is therefore worth using an
elimination method.

\begin{figure}
\begin{center}
\includegraphics[width=0.95\linewidth]{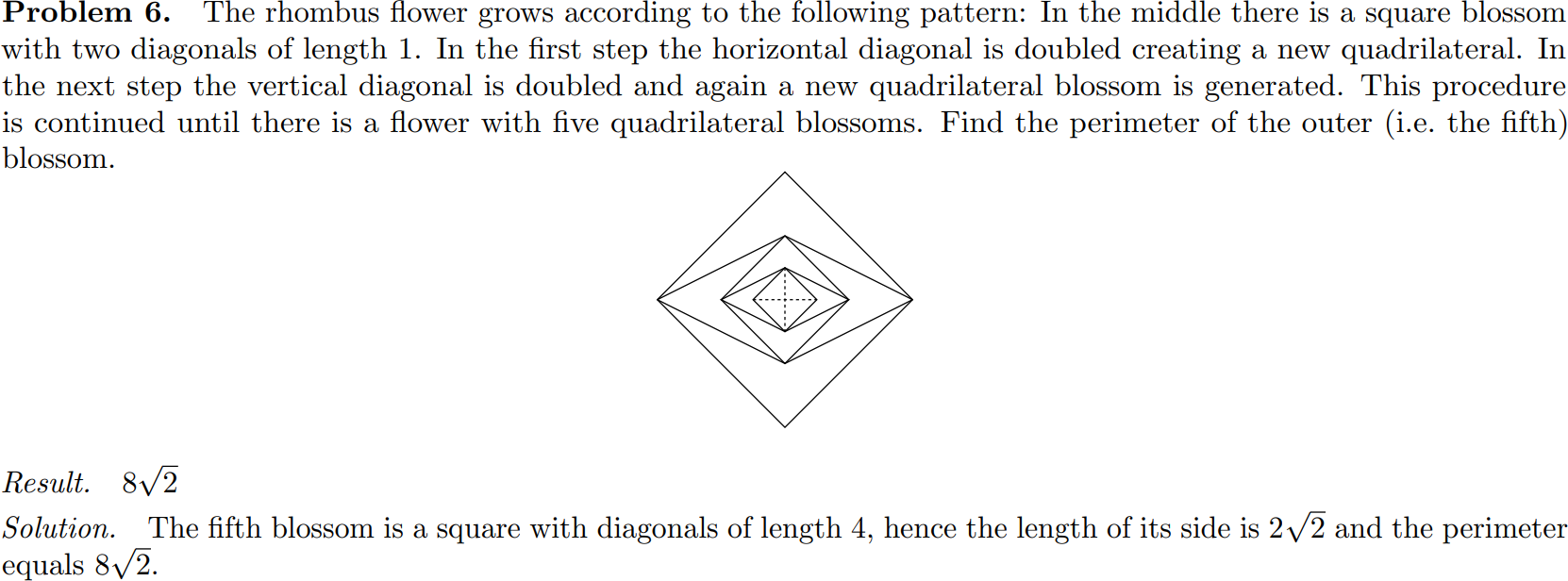}
\caption{Problem setting 6 and the official solution of Náboj 2023}\label{Problem6} 
\end{center}
\end{figure}

As an illustration, we show how the program solves Problem 6 of Náboj 2023 (Fig.~\ref{Problem6}).

\begin{figure}[h]
\begin{center}
\includegraphics[width=\linewidth]{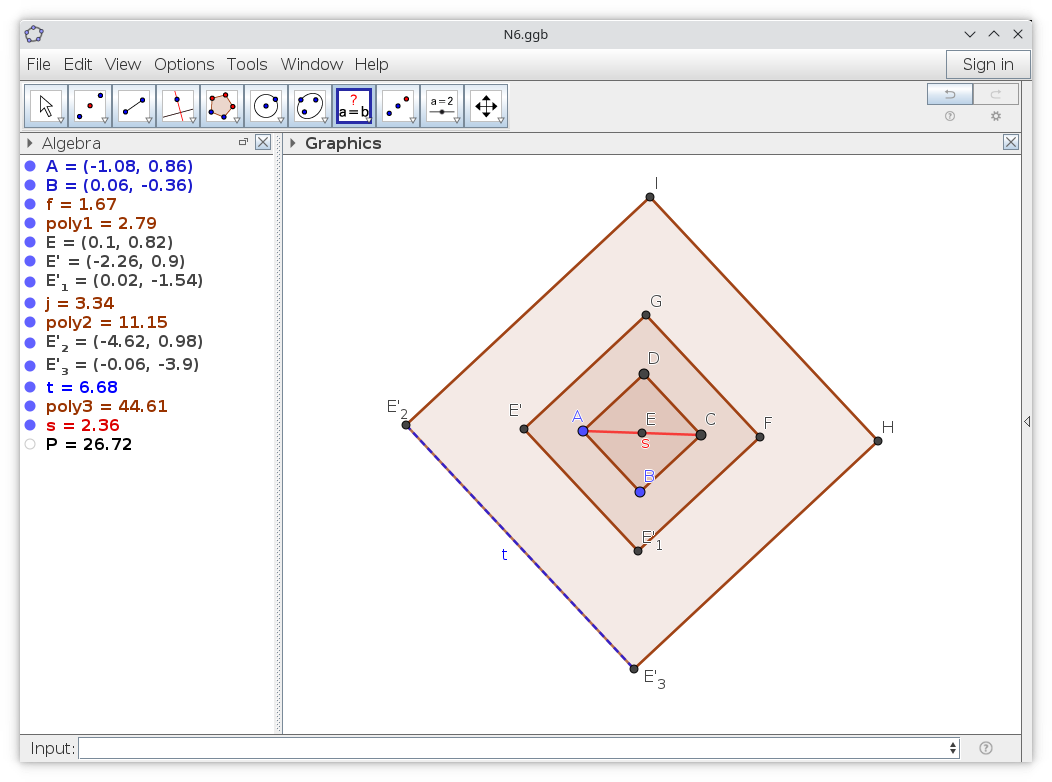}
\caption{Sketching Problem 6 in GeoGebra Discovery}\label{GD6} 
\end{center}
\end{figure}

\begin{figure}
\begin{center}
\includegraphics[width=0.5\linewidth]{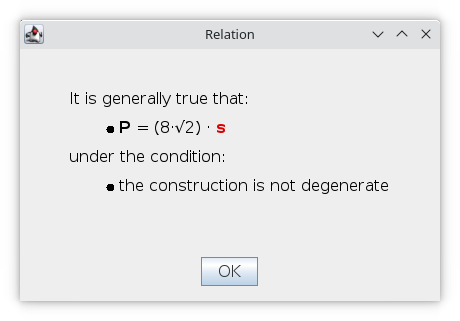}
\caption{Report of a symbolic analysis of Problem 6 in GeoGebra Discovery}\label{GD6r} 
\end{center}
\end{figure}

When drawing the figure in GeoGebra Discovery (Fig.~\ref{GD6}), we learn that the problem can be simplified to
three squares. First, an arbitrary square $ABCD$ is drawn. Then midpoint $E$ of $AC$ is defined.
This point will be then reflected about $A$ to get $E'$ and about $B$ to get $E_1'$. Now a second
square $E'E_1'FG$ is drawn. Finally, by reflecting $E$ about $E'$ and $E_1'$ we get points
$E_2'$ and $E_3'$, respectively, and create the square $E_2'E_3'HI$ as well. By defining $s=AC$,
$t=E_2'E_3'$ and $P=4t$, we can use GeoGebra Discovery's Relation tool to compare $s$ and $P$ and
we learn (after pressing the button ``More$\ldots$'' to obtain a symbolic analysis) that
$P=8\sqrt2\cdot s$. The report of the symbolic analysis of Problem 6 in
GeoGebra Discovery shows that ``It is generally true that:
$P=(8\sqrt{2})\cdot s$ under the condition: the construction is not
degenerate'' (see Fig.~\ref{GD6r}). The construction steps in GeoGebra
Discovery for Problem 6 can be taken from Table \ref{GD6cp}. (These extra
explanations are listed in the Appendix.)

\definecolor{ggbcolor0}{RGB}{27, 27, 205}
\definecolor{ggbcolor1}{RGB}{27, 27, 205}
\definecolor{ggbcolor2}{RGB}{153, 51, 0}
\definecolor{ggbcolor3}{RGB}{153, 51, 0}
\definecolor{ggbcolor4}{RGB}{153, 51, 0}
\definecolor{ggbcolor5}{RGB}{68, 68, 68}
\definecolor{ggbcolor6}{RGB}{68, 68, 68}
\definecolor{ggbcolor7}{RGB}{153, 51, 0}
\definecolor{ggbcolor8}{RGB}{153, 51, 0}
\definecolor{ggbcolor9}{RGB}{68, 68, 68}
\definecolor{ggbcolor10}{RGB}{68, 68, 68}
\definecolor{ggbcolor11}{RGB}{68, 68, 68}
\definecolor{ggbcolor12}{RGB}{153, 51, 0}
\definecolor{ggbcolor13}{RGB}{153, 51, 0}
\definecolor{ggbcolor14}{RGB}{153, 51, 0}
\definecolor{ggbcolor15}{RGB}{68, 68, 68}
\definecolor{ggbcolor16}{RGB}{68, 68, 68}
\definecolor{ggbcolor17}{RGB}{153, 51, 0}
\definecolor{ggbcolor18}{RGB}{153, 51, 0}
\definecolor{ggbcolor19}{RGB}{68, 68, 68}
\definecolor{ggbcolor20}{RGB}{68, 68, 68}
\definecolor{ggbcolor21}{RGB}{153, 51, 0}
\definecolor{ggbcolor22}{RGB}{0, 0, 255}
\definecolor{ggbcolor23}{RGB}{153, 51, 0}
\definecolor{ggbcolor24}{RGB}{68, 68, 68}
\definecolor{ggbcolor25}{RGB}{68, 68, 68}
\definecolor{ggbcolor26}{RGB}{153, 51, 0}
\definecolor{ggbcolor27}{RGB}{153, 51, 0}
\definecolor{ggbcolor28}{RGB}{220, 0, 0}
Here we highlight that the problem setting could be further simplified by skipping the construction
of the two latter squares. In fact, only the reflection points matter. Also, we used the fact
that the perimeter of a square equals to four times the length of a side, but this piece of
information could have been ignored and asked the program to learn this on its own.


At this point, we will jump to the last steps and assuming that it is possible to do that without loss of generality, GeoGebra Discovery, in the background (in an invisible way for the normal user, but in a verifiable way via its debug messages), decides to substitute $A = (0,0)$.
This will simplify the computations and the final question is if $e_{21}:m\cdot v_{24}=4v_{23}$ holds
for a given value of $m$. This can be answered by eliminating all variables but $m$ from all of the equations $e_1, e_2, \ldots, e_{14}, e_{18}, \ldots, e_{21}$, and we learn that $m^2=128$.
That is, by assuming that $m>0$, we obtain that $m=8\sqrt2$.

The algebraic solution is, of course, quite complicated. Also, the way used in GeoGebra Discovery
by constructing all the required objects, may still be very complicated when compared to the
quick official solution. 

\section{Problems that can be Solved with GeoGebra Discovery}

In this section we list four additional problems that can be solved with GeoGebra Discovery,
assuming some effort. In fact, some other Náboj 2023 problems can be supported as well,
but they may require some additional steps. See the next section for more details.

\subsection{Pentominos (Problem 15)}

\begin{figure}
\begin{center}
\includegraphics[width=0.95\linewidth]{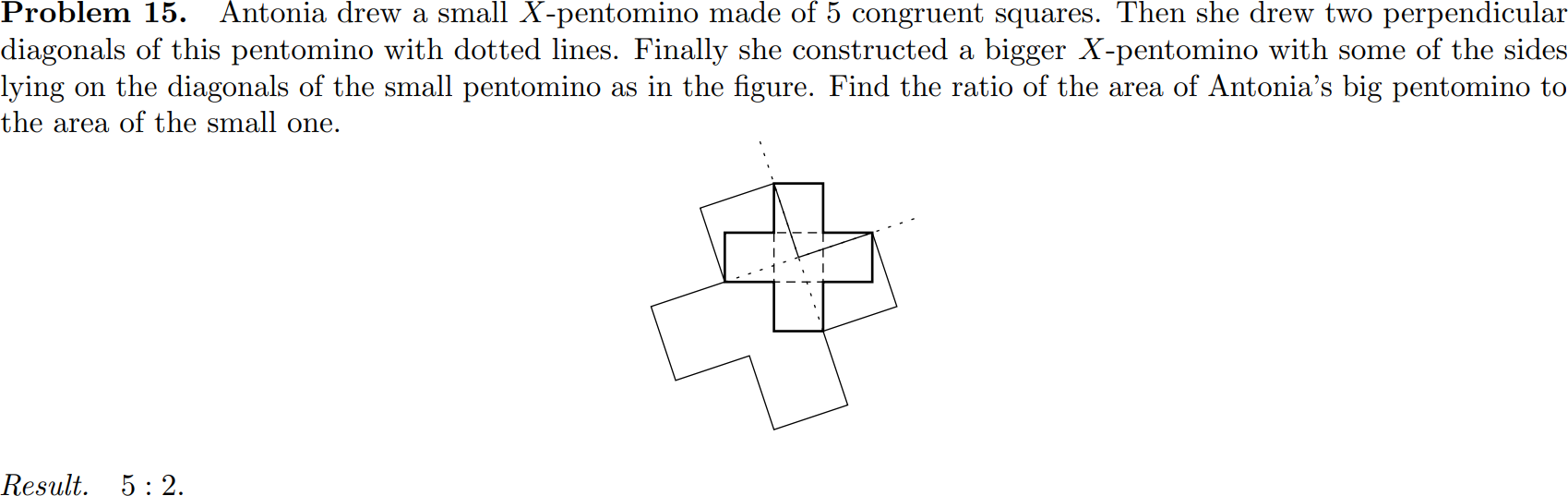}\\
\caption{Problem setting 15 
}\label{Problem15} 
\end{center}
\end{figure}

The problem setting can be seen in Fig. \ref{Problem15}.

\begin{figure}
\begin{center}
\includegraphics[width=0.5\linewidth]{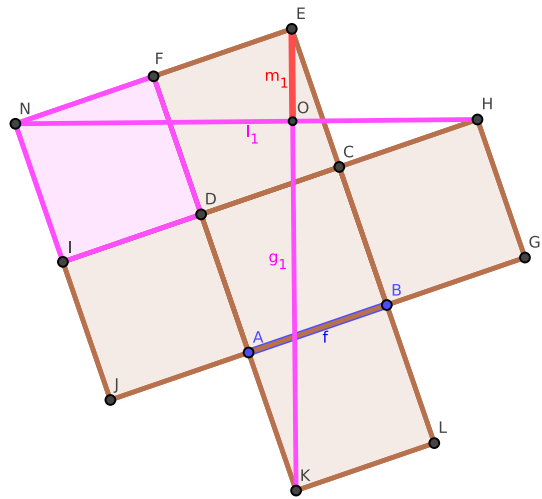}
\caption{Problem setting 15 in GeoGebra Discovery}\label{GD15} 
\end{center}
\end{figure}

At a first look, it seems complicated to draw a figure that describes the problem setting adequately.
Some attempts may lead to Fig. \ref{GD15}: lines $g_1=EK$ and (after extending the large
pentomine with square $DINF$) $l_1=HN$ help finding point $O$.
Then, segment $m_1=EO$ will be one side of the small $X$-pentomino, and it will be possible
to compare it to one of the sides of the large pentomino. For more details about the construction steps have a look at Table \ref{GD15cp}.

\definecolor{ggbcolor0}{RGB}{27, 27, 205}
\definecolor{ggbcolor1}{RGB}{27, 27, 205}
\definecolor{ggbcolor2}{RGB}{153, 51, 0}
\definecolor{ggbcolor3}{RGB}{0, 0, 255}
\definecolor{ggbcolor4}{RGB}{153, 51, 0}
\definecolor{ggbcolor5}{RGB}{68, 68, 68}
\definecolor{ggbcolor6}{RGB}{68, 68, 68}
\definecolor{ggbcolor7}{RGB}{153, 51, 0}
\definecolor{ggbcolor8}{RGB}{153, 51, 0}
\definecolor{ggbcolor9}{RGB}{153, 51, 0}
\definecolor{ggbcolor10}{RGB}{153, 51, 0}
\definecolor{ggbcolor11}{RGB}{153, 51, 0}
\definecolor{ggbcolor12}{RGB}{68, 68, 68}
\definecolor{ggbcolor13}{RGB}{68, 68, 68}
\definecolor{ggbcolor14}{RGB}{153, 51, 0}
\definecolor{ggbcolor15}{RGB}{153, 51, 0}
\definecolor{ggbcolor16}{RGB}{153, 51, 0}
\definecolor{ggbcolor17}{RGB}{153, 51, 0}
\definecolor{ggbcolor18}{RGB}{153, 51, 0}
\definecolor{ggbcolor19}{RGB}{68, 68, 68}
\definecolor{ggbcolor20}{RGB}{68, 68, 68}
\definecolor{ggbcolor21}{RGB}{153, 51, 0}
\definecolor{ggbcolor22}{RGB}{153, 51, 0}
\definecolor{ggbcolor23}{RGB}{153, 51, 0}
\definecolor{ggbcolor24}{RGB}{153, 51, 0}
\definecolor{ggbcolor25}{RGB}{153, 51, 0}
\definecolor{ggbcolor26}{RGB}{68, 68, 68}
\definecolor{ggbcolor27}{RGB}{68, 68, 68}
\definecolor{ggbcolor28}{RGB}{153, 51, 0}
\definecolor{ggbcolor29}{RGB}{153, 51, 0}
\definecolor{ggbcolor30}{RGB}{153, 51, 0}
\definecolor{ggbcolor31}{RGB}{153, 51, 0}
\definecolor{ggbcolor32}{RGB}{153, 51, 0}
\definecolor{ggbcolor33}{RGB}{68, 68, 68}
\definecolor{ggbcolor34}{RGB}{68, 68, 68}
\definecolor{ggbcolor35}{RGB}{153, 51, 0}
\definecolor{ggbcolor36}{RGB}{153, 51, 0}
\definecolor{ggbcolor37}{RGB}{130, 0, 130}
\definecolor{ggbcolor38}{RGB}{130, 0, 130}
\definecolor{ggbcolor39}{RGB}{130, 0, 130}
\definecolor{ggbcolor40}{RGB}{130, 0, 130}
\definecolor{ggbcolor41}{RGB}{68, 68, 68}
\definecolor{ggbcolor42}{RGB}{68, 68, 68}
\definecolor{ggbcolor43}{RGB}{130, 0, 130}
\definecolor{ggbcolor44}{RGB}{130, 0, 130}
\definecolor{ggbcolor45}{RGB}{130, 0, 130}
\definecolor{ggbcolor46}{RGB}{68, 68, 68}
\definecolor{ggbcolor47}{RGB}{220, 0, 0}

Now, asking the relation between $m_1$ and $f$ we obtain that $f=\frac12\cdot\sqrt{10}\cdot m_1$.
Even if GeoGebra Discovery cannot compare areas symbolically in a direct way, we can still conclude
that the ratio of the areas must be
$$f^2:{m_1}^2=10:4=5:2.$$

We remark here that a sophisticated way to do the construction may give a quicker result
than the official solution. It may be, however, not trivial to find this alternative solution.

\subsection{A Right Triangle (Problem 23)}
\begin{figure}
\begin{center}
\includegraphics[width=0.95\linewidth]{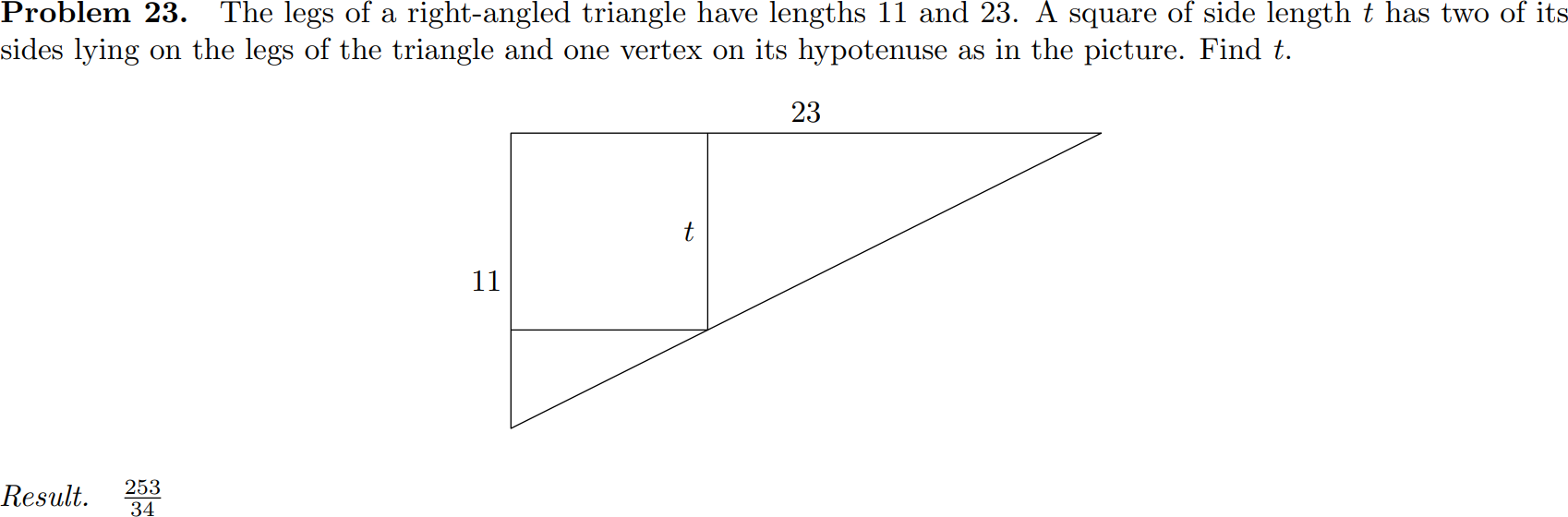}
\caption{Problem setting 23 
}\label{Problem23} 
\end{center}
\end{figure}

The problem setting can be seen in Fig. \ref{Problem23}.


\begin{figure}
\begin{center}
\includegraphics[width=0.7\linewidth]{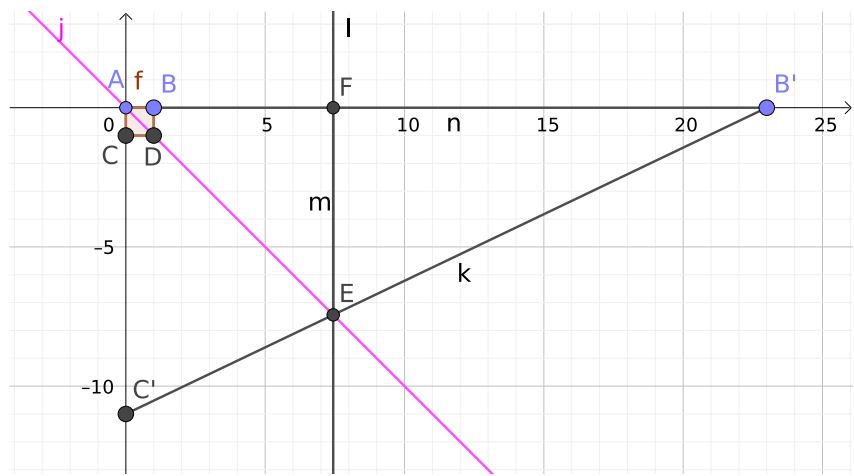}\\
\caption{Sketch for Problem 23}\label{GD23} 
\end{center}
\end{figure}

We use some recent features of GeoGebra Discovery to solve this problem.
Most importantly, a square is created based on free points $A=(0,0)$ and $B=(1,0)$.
They must not be defined with the help of any axes, because in that case the background
proof will fail and no output will be obtained. Here, instead of copying the initial
square several times, we use the Dilate tool to stretch the segment $AB$ and $AC$
to get $B'$ and $C'$ accordingly. Another trick is to create the diagonal $j=AD$ of the initial square.
Now the intersection $E$ of $k=B'C'$ and $j$ is the searched point. Finally,
projecting $E$ on $n=AB'$ and obtaining intersection point $F$ of perpendicular $l$ and line $n$,
comparison of $m=EF$ and $f=AB$ is to be done. And, indeed $m=253/34$, as expected.
(See the sketch in Fig.~\ref{GD23} in GeoGebra Discovery.)

\definecolor{ggbcolor0}{RGB}{40, 40, 170}
\definecolor{ggbcolor1}{RGB}{40, 40, 170}
\definecolor{ggbcolor2}{RGB}{40, 40, 170}
\definecolor{ggbcolor4}{RGB}{153, 51, 0}
\definecolor{ggbcolor5}{RGB}{153, 51, 0}
\definecolor{ggbcolor6}{RGB}{153, 51, 0}
\definecolor{ggbcolor7}{RGB}{68, 68, 68}
\definecolor{ggbcolor8}{RGB}{68, 68, 68}
\definecolor{ggbcolor9}{RGB}{153, 51, 0}
\definecolor{ggbcolor10}{RGB}{153, 51, 0}
\definecolor{ggbcolor11}{RGB}{130, 0, 130}
\definecolor{ggbcolor12}{RGB}{68, 68, 68}
\definecolor{ggbcolor14}{RGB}{68, 68, 68}
\definecolor{ggbcolor16}{RGB}{68, 68, 68}
\definecolor{ggbcolor17}{RGB}{220, 0, 0}

Here we remark that finding the rational value $253/34$ (it is approximately 7.44) seems very difficult
unless one does not solve the problem explicitly (as shown in the official solution, by using an equation).
If a user has some routine in GeoGebra, sketching the problem may take a shorter time than finding
the required equation (even if it is a linear one). The construction steps can be found in Table \ref{GD23cp}.

\subsection{A Triangle and a Circle (Problem 47)}

Problem 47 was not even accessible during the contest for most of the teams because it was almost the last
problem in the list and they were not that fast.

\begin{figure}
\begin{center}
\includegraphics[width=0.95\linewidth]{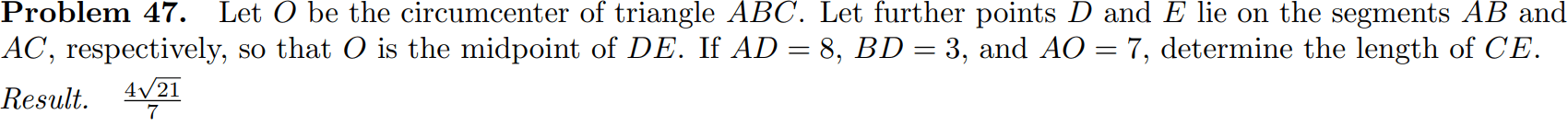}\\
\caption{Problem setting 47
}\label{Problem47} 
\end{center}
\end{figure}

The problem setting can be seen in Fig. \ref{Problem47}.

The official solution of Problem 47 required some non-trivial ideas.
When using GeoGebra Discovery, we may face the question how the problem setting can be constructed, which is shown on Table \ref{GD47cp}. Since the lengths $AD=8$ and $BD=3$ are given, it seems reasonable to create $AB$ arbitrarily
and create $D$ by using the command \texttt{Dilate(B, 8/11, A)}. Now, we create another
point $A'$ in a similar way, by placing $A'$ on $AB$ and letting $AA'=7$. This helps us
restricting the position of $O$ because it must be on the circumcircle of the circle $c$ with center $A$
and radius $AA'$. On the other hand, $O$ must lie on the perpendicular bisector $g$ of $AB$.
At this point we already know the position of $O=c\cap g$. (In fact, there may be two solutions here,
but they are identical in the sense of symmetry.)

Now, by reflecting $D$ about $O$ we obtain $E$. By having $E$, we already know the line $BC$.
To get the point $C$ we only have to intersect this line with the circumcircle $c$. (Again, there
are two solutions, but the other one $C'$ leads to a degenerate case because it yields $A=C'$.
To force getting the non-degenerate case we need to click near the intersection point with the mouse.
Otherwise GeoGebra Discovery will compute with \textit{both} cases at the same time.)

As a final step, we designate the unit length. Luckily, $DA'$ is exactly $1$. So we just have to compare $j=AE$ and $i=DA'$. As expected, the result is $j=4/7\cdot\sqrt{21}\cdot i$. Thus, $CE=j=\frac{4\sqrt{21}}{7}$.

\begin{figure}
\begin{center}
\includegraphics[width=0.5\linewidth]{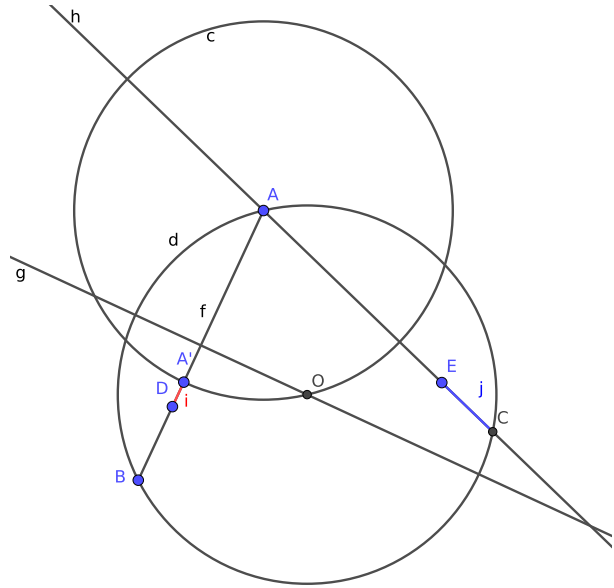}\\
\caption{Sketch for Problem 47}\label{GD47} 
\end{center}
\end{figure}

The sketch can be seen in Fig. \ref{GD47}.

\definecolor{ggbcolor0}{RGB}{27, 27, 205}
\definecolor{ggbcolor1}{RGB}{27, 27, 205}
\definecolor{ggbcolor3}{RGB}{27, 27, 205}
\definecolor{ggbcolor4}{RGB}{27, 27, 205}
\definecolor{ggbcolor6}{RGB}{68, 68, 68}
\definecolor{ggbcolor8}{RGB}{27, 27, 205}
\definecolor{ggbcolor11}{RGB}{220, 0, 0}
\definecolor{ggbcolor12}{RGB}{68, 68, 68}
\definecolor{ggbcolor13}{RGB}{0, 0, 255}

\section{Problems that Require Further Improvements}

In this section we take an overview of other examples that cannot be fully solved
in simple steps in Geo\-Gebra Discovery. Some hints may be, however, obtained.
Instead of getting such hints, we summarize how the software tool could be extended
to be able to give full solutions for such problems.

\begin{figure}
\begin{center}
\includegraphics[width=0.95\linewidth]{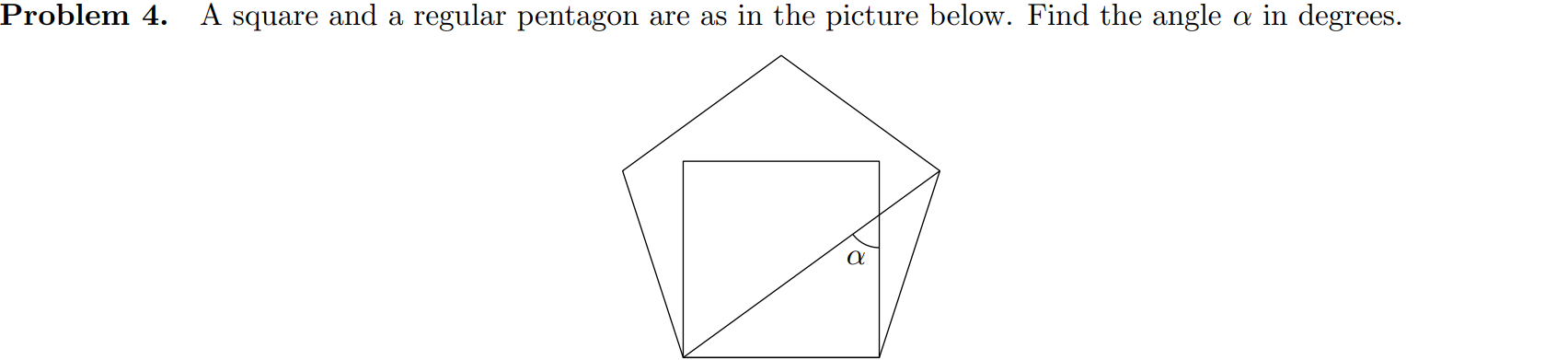}\\
\caption{Problem 4}\label{Problem4}
\end{center}
\end{figure}

\begin{figure}
\begin{center}
\includegraphics[width=0.95\linewidth]{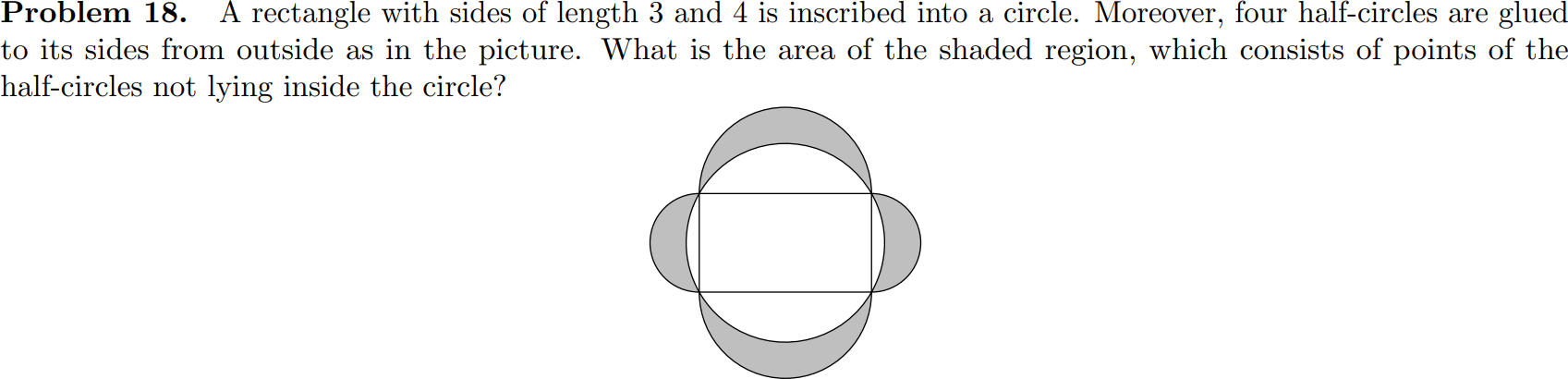}\\
\caption{Problem 18}\label{Problem18}
\end{center}
\end{figure}

\begin{figure}
\begin{center}
\includegraphics[width=0.95\linewidth]{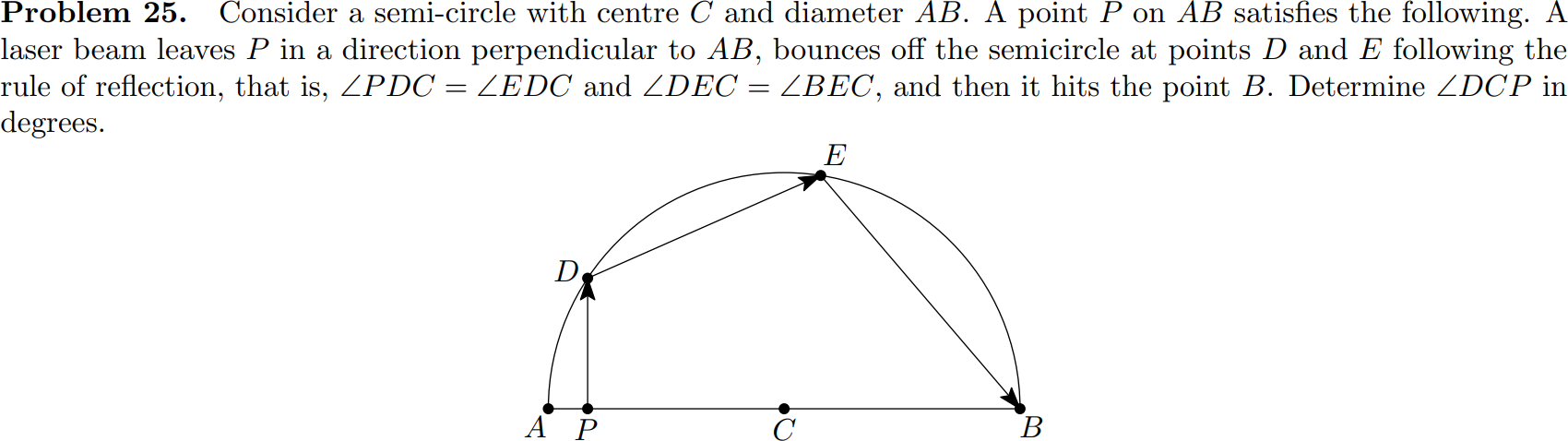}\\
\caption{Problem 25}\label{Problem25}
\end{center}
\end{figure}

\begin{figure}
\begin{center}
\includegraphics[width=0.95\linewidth]{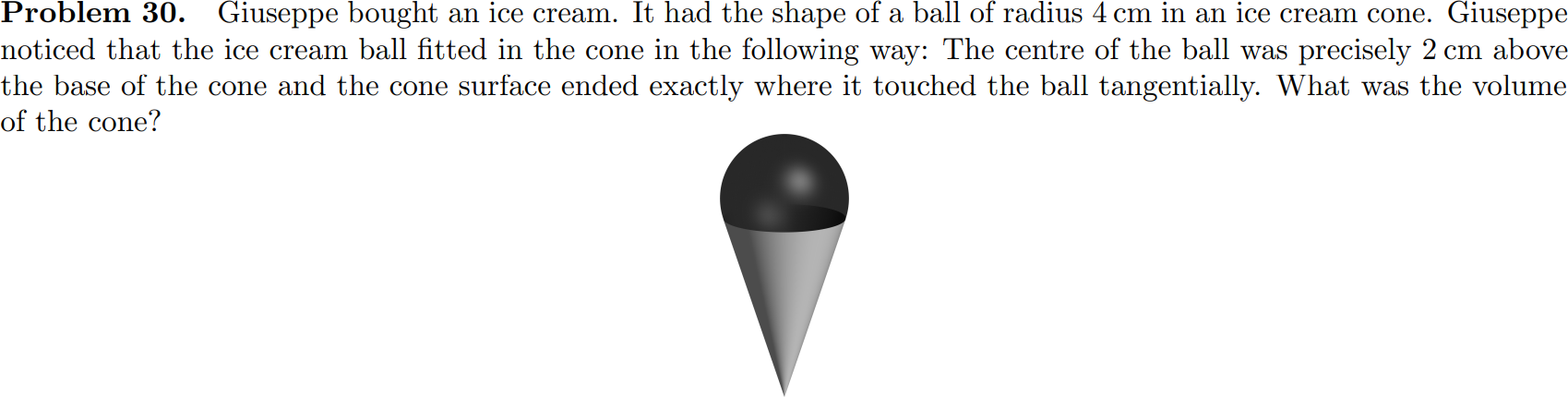}\\
\caption{Problem 30}\label{Problem30}
\end{center}
\end{figure}

\begin{figure}
\begin{center}
\includegraphics[width=0.95\linewidth]{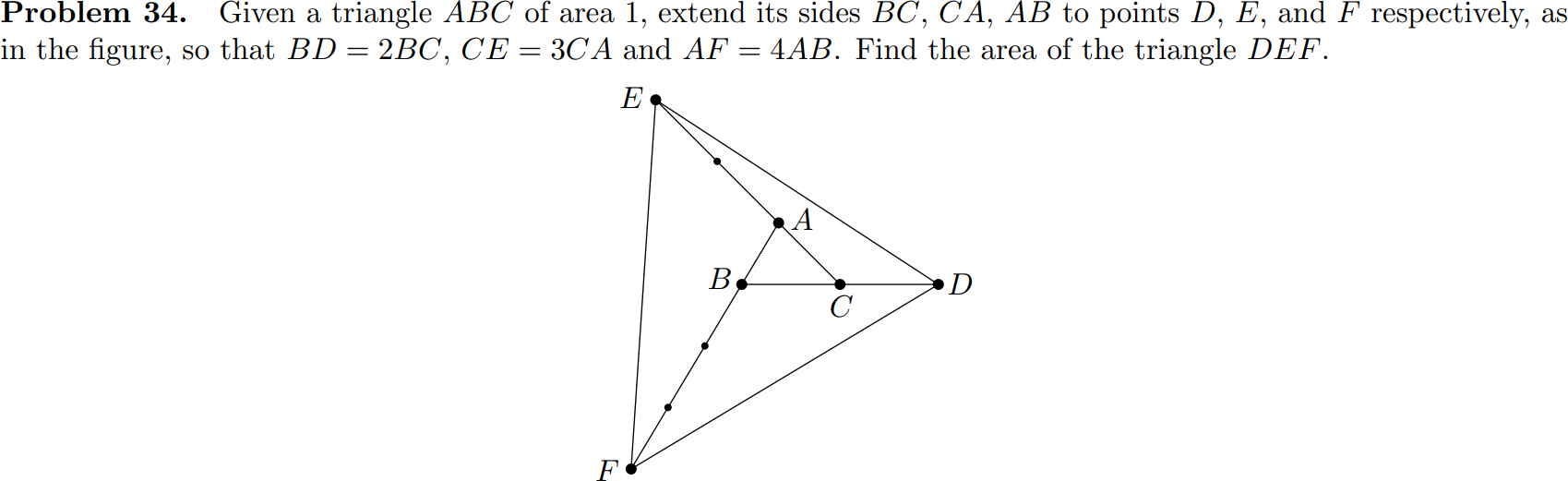}\\
\caption{Problem 34}\label{Problem34}
\end{center}
\end{figure}

\begin{figure}
\begin{center}
\includegraphics[width=0.95\linewidth]{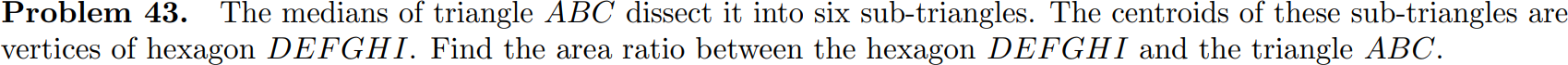}\\
\caption{Problem 43}\label{Problem43}
\end{center}
\end{figure}

\begin{figure}
\begin{center}
\includegraphics[width=0.95\linewidth]{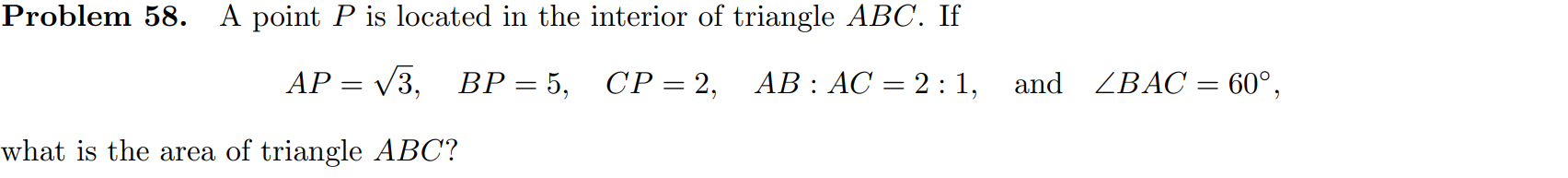}\\
\caption{Problem 58}\label{Problem58}
\end{center}
\end{figure}

A first set of problems (4, 25 and 58, see Fig. \ref{Problem4}, \ref{Problem25} and \ref{Problem58})
are related with angles, another set (Problems
18, 30, 34, 43 and 58, see Fig. \ref{Problem18}, \ref{Problem30},
\ref{Problem34}, \ref{Problem43} and \ref{Problem58}) is about areas. Angle support (via symbolic computation)
is very poor in GeoGebra Discovery: this has roots in a non-bijective relationship
between angles and their algebraic counterparts. Area support is also somewhat minimal,
because it is restricted to triangles, and the expected way of use is not polished yet
in the software tool.

One can find that Problems 25 and 58 have some common roots. 
They can be formulated with ``implicit assumptions''.

\begin{figure}
\begin{center}
\includegraphics[width=0.7\linewidth]{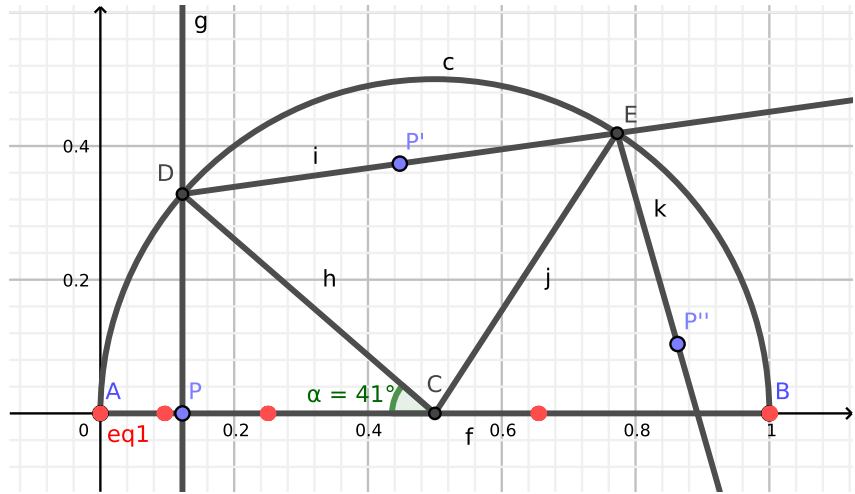}\\
\caption{A possible approach to solve Problem 25}\label{GD25}
\end{center}
\end{figure}

We have a closer look at Problem 25. Let point $F$ be the next bounce, when we assume that a laser beam starts
from point $P$. Now, a command like \texttt{LocusEquation($F$==$B$,$\angle DCP$)} could
address the question (but having angles in the second parameter is not implemented).
In fact, we may already get the exact position of $P$
when applying consecutive reflections. According to Fig. \ref{GD25}, when
reflecting $P$ about $CD$, and intersecting the line connecting $D$ and the mirror image $P'$
with the semicircle, we can obtain $E$. Another reflection can yield $P''$ and the final
visualization can be achieved with \texttt{LocusEquation(AreCollinear($E$,$P''$,$B$),$P$)}.
Since GeoGebra Discovery shows 5 isolated points, one can conjecture that there is something
to do with a regular pentagon. Here, unfortunately, the factorization of the obtained
polynomial does not help, because the interesting quadratic numbers are appearing
just approximately. A deeper symbolic study shows that (by assuming $A=(0,0)$ and $B=(1,0)$)
for the $x$-coordinate of $P$ is one of the roots of the polynomial
$64x^5-128x^4+80x^3-17x^2+x$, and they are
$$0,\frac{1}{16}\cdot(-2\cdot\sqrt{5}+6),\frac14,\frac{1}{16}\cdot(2\cdot\sqrt{5}+6),1,$$
and to these values belong the $\alpha$ values
$$0^{\textrm{o}},36^{\textrm{o}},60^{\textrm{o}},288^{\textrm{o}},360^{\textrm{o}},$$
the latter two ones without real geometrical meaning. Finally we can conclude that $\alpha=36^{\textrm{o}}$,
this is the only meaningful solution. But, all of this derivation requires
some additional steps, GeoGebra Discovery alone does not bring a satisfactory final answer.

\begin{figure}
\begin{center}
\includegraphics[width=0.9\linewidth]{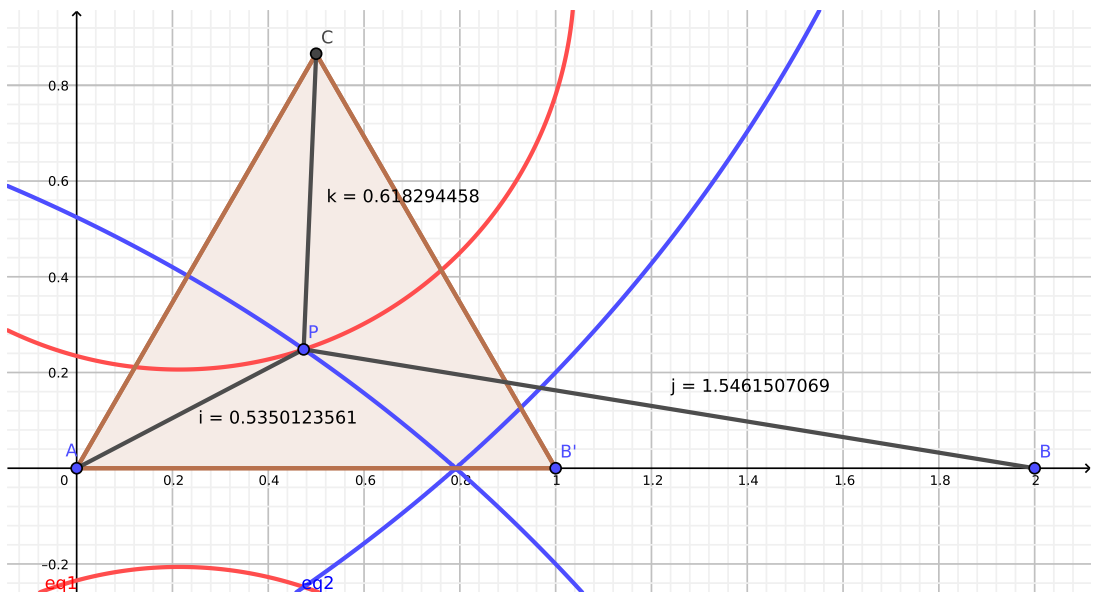}\\
\caption{Obtaining an incorrect conjecture to solve Problem 58}\label{GD58}
\end{center}
\end{figure}

Finally, we show a \textit{wrong} conjecture for Problem 58, based on GeoGebra Discovery.
Like Problem 25, an implicit locus equation seems here helpful. Let us, first,
create a regular triangle $AB'C$ with $A=(0,0)$, $B=(1,0)$, and reflect $A$ about $B'$ to get $B$. Clearly,
these preparations are sufficient to ensure assumptions $AB:AC=2:1$ and $\angle BAC=60^{\textrm{o}}$
(See Fig. \ref{GD58}).
This dummy triangle has the area $\frac{\sqrt3}2$.
Now, we create an arbitrary point $P$ and connect it with points $A$, $B$ and $C$, to get segments
$i$, $j$ and $k$, respectively. We create two locus equations with the commands
\texttt{LocusEquation($j/k$==5/2,$P$)} and
\texttt{LocusEquation($i/k$==sqrt(3)/2,$P$)}. Now, we want to find the correct position for $P$,
so we consider the intersection of the two locus curves visually. After zooming in, we learn that
for $P=(0.4739140532,0.24828147621)$ we obtain $k=0.618294458$ which seems to be close enough to the well known
number $f=\frac{\sqrt5-1}2$. If so, the triangle must be enlarged by a factor $1/f\cdot2$ which is twice the golden ratio, $2\varphi=\sqrt5+1$. Finally,
the triangle will have the area $\frac{\sqrt3}2\cdot(\sqrt5+1)\approx9.06913\ldots$

Assuming that the golden ratio plays a role here is, however, incorrect. The correct solution is $\frac{6+7\cdot\sqrt3}2\approx9.06217\ldots$\footnote{%
See \url{https://math.old.naboj.org/archive/problems/pdf/math/2023_en_sol.pdf} for a full computation.}
Note that the ratio between the two values is $1.000767\ldots$ which is remarkably small difference.

\section{Conclusion}

The last section showed that GeoGebra Discovery can be a useful tool to get a correct conjecture if
the right steps are taken to finish the solution, but it can also be misleading in some delicate situations.
On the other hand, several contest problems can be handled and solved with minimal effort by using this tool.
We need to admit that a good knowledge of the software is unavoidable. However, experienced users may
need just a couple of steps to achieve the solution.

For a future improvement, full support of computing angles and areas seems to be a great step forward.
Some sophisticated problems, however, may need further developments towards symbolic computations that are
based on implicit assumptions. 

\section{Acknowledgements}
The second author was partially
supported by a grant
PID2020-113192GB-I00 (Mathematical Visualization: Foundations, Algorithms
and Applications) from the Spanish MICINN.

\bibliography{kovzol,external}

\section{Appendix}
In this Appendix, we provide the construction protocols in GeoGebra Discovery for the problems presented above.

\definecolor{ggbcolor0}{RGB}{27, 27, 205}
\definecolor{ggbcolor1}{RGB}{27, 27, 205}
\definecolor{ggbcolor2}{RGB}{153, 51, 0}
\definecolor{ggbcolor3}{RGB}{153, 51, 0}
\definecolor{ggbcolor4}{RGB}{153, 51, 0}
\definecolor{ggbcolor5}{RGB}{68, 68, 68}
\definecolor{ggbcolor6}{RGB}{68, 68, 68}
\definecolor{ggbcolor7}{RGB}{153, 51, 0}
\definecolor{ggbcolor8}{RGB}{153, 51, 0}
\definecolor{ggbcolor9}{RGB}{68, 68, 68}
\definecolor{ggbcolor10}{RGB}{68, 68, 68}
\definecolor{ggbcolor11}{RGB}{68, 68, 68}
\definecolor{ggbcolor12}{RGB}{153, 51, 0}
\definecolor{ggbcolor13}{RGB}{153, 51, 0}
\definecolor{ggbcolor14}{RGB}{153, 51, 0}
\definecolor{ggbcolor15}{RGB}{68, 68, 68}
\definecolor{ggbcolor16}{RGB}{68, 68, 68}
\definecolor{ggbcolor17}{RGB}{153, 51, 0}
\definecolor{ggbcolor18}{RGB}{153, 51, 0}
\definecolor{ggbcolor19}{RGB}{68, 68, 68}
\definecolor{ggbcolor20}{RGB}{68, 68, 68}
\definecolor{ggbcolor21}{RGB}{153, 51, 0}
\definecolor{ggbcolor22}{RGB}{0, 0, 255}
\definecolor{ggbcolor23}{RGB}{153, 51, 0}
\definecolor{ggbcolor24}{RGB}{68, 68, 68}
\definecolor{ggbcolor25}{RGB}{68, 68, 68}
\definecolor{ggbcolor26}{RGB}{153, 51, 0}
\definecolor{ggbcolor27}{RGB}{153, 51, 0}
\definecolor{ggbcolor28}{RGB}{220, 0, 0}

\begin{longtable}[]{rlcp{4cm}}
\caption{Construction protocol in GeoGebra Discovery for Problem~6}\label{GD6cp}\\
\toprule
No.&Name&Toolbar Icon&Description\tabularnewline
\midrule
\endhead
\textcolor{ggbcolor2}{1}&\textcolor{ggbcolor2}{Polygon poly1}&\textcolor{ggbcolor2}{\raisebox{-1mm}{\includegraphics[width=0.5cm]{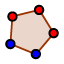}}}&\textcolor{ggbcolor2}{Polygon(A, B, 4)}\tabularnewline
\textcolor{ggbcolor9}{2}&\textcolor{ggbcolor9}{Point E}&\textcolor{ggbcolor9}{\raisebox{-1mm}{\includegraphics[width=0.5cm]{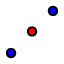}}}&\textcolor{ggbcolor9}{Midpoint of A, C}\tabularnewline
\textcolor{ggbcolor10}{3}&\textcolor{ggbcolor10}{Point E'}&\textcolor{ggbcolor10}{\raisebox{-1mm}{\includegraphics[width=0.5cm]{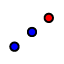}}}&\textcolor{ggbcolor10}{E mirrored at A}\tabularnewline
\textcolor{ggbcolor11}{4}&\textcolor{ggbcolor11}{Point E'\_{1}}&\textcolor{ggbcolor11}{\raisebox{-1mm}{\includegraphics[width=0.5cm]{mode_mirroratpoint}}}&\textcolor{ggbcolor11}{E mirrored at B}\tabularnewline
\textcolor{ggbcolor12}{5}&\textcolor{ggbcolor12}{Polygon poly2}&\textcolor{ggbcolor12}{\raisebox{-1mm}{\includegraphics[width=0.5cm]{mode_regularpolygon}}}&\textcolor{ggbcolor12}{Polygon(E', E'\_{1}, 4)}\tabularnewline
\textcolor{ggbcolor19}{6}&\textcolor{ggbcolor19}{Point E'\_{2}}&\textcolor{ggbcolor19}{\raisebox{-1mm}{\includegraphics[width=0.5cm]{mode_mirroratpoint}}}&\textcolor{ggbcolor19}{E mirrored at E'}\tabularnewline
\textcolor{ggbcolor20}{7}&\textcolor{ggbcolor20}{Point E'\_{3}}&\textcolor{ggbcolor20}{\raisebox{-1mm}{\includegraphics[width=0.5cm]{mode_mirroratpoint}}}&\textcolor{ggbcolor20}{E mirrored at E'\_{1}}\tabularnewline
\textcolor{ggbcolor21}{8}&\textcolor{ggbcolor21}{Polygon poly3}&\textcolor{ggbcolor21}{\raisebox{-1mm}{\includegraphics[width=0.5cm]{mode_regularpolygon}}}&\textcolor{ggbcolor21}{Polygon(E'\_{2}, E'\_{3}, 4)}\tabularnewline
\textcolor{ggbcolor22}{9}&\textcolor{ggbcolor22}{Segment t}&\textcolor{ggbcolor22}{\raisebox{-1mm}{\includegraphics[width=0.5cm]{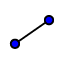}}}&\textcolor{ggbcolor22}{Segment E'\_{2}, E'\_{3}}\tabularnewline
\textcolor{ggbcolor28}{10}&\textcolor{ggbcolor28}{Segment s}&\textcolor{ggbcolor28}{\raisebox{-1mm}{\includegraphics[width=0.5cm]{mode_segment}}}&\textcolor{ggbcolor28}{Segment A, C}\tabularnewline
11&Number P& &4t\tabularnewline
\bottomrule
\end{longtable}

\definecolor{ggbcolor0}{RGB}{27, 27, 205}
\definecolor{ggbcolor1}{RGB}{27, 27, 205}
\definecolor{ggbcolor2}{RGB}{153, 51, 0}
\definecolor{ggbcolor3}{RGB}{0, 0, 255}
\definecolor{ggbcolor4}{RGB}{153, 51, 0}
\definecolor{ggbcolor5}{RGB}{68, 68, 68}
\definecolor{ggbcolor6}{RGB}{68, 68, 68}
\definecolor{ggbcolor7}{RGB}{153, 51, 0}
\definecolor{ggbcolor8}{RGB}{153, 51, 0}
\definecolor{ggbcolor9}{RGB}{153, 51, 0}
\definecolor{ggbcolor10}{RGB}{153, 51, 0}
\definecolor{ggbcolor11}{RGB}{153, 51, 0}
\definecolor{ggbcolor12}{RGB}{68, 68, 68}
\definecolor{ggbcolor13}{RGB}{68, 68, 68}
\definecolor{ggbcolor14}{RGB}{153, 51, 0}
\definecolor{ggbcolor15}{RGB}{153, 51, 0}
\definecolor{ggbcolor16}{RGB}{153, 51, 0}
\definecolor{ggbcolor17}{RGB}{153, 51, 0}
\definecolor{ggbcolor18}{RGB}{153, 51, 0}
\definecolor{ggbcolor19}{RGB}{68, 68, 68}
\definecolor{ggbcolor20}{RGB}{68, 68, 68}
\definecolor{ggbcolor21}{RGB}{153, 51, 0}
\definecolor{ggbcolor22}{RGB}{153, 51, 0}
\definecolor{ggbcolor23}{RGB}{153, 51, 0}
\definecolor{ggbcolor24}{RGB}{153, 51, 0}
\definecolor{ggbcolor25}{RGB}{153, 51, 0}
\definecolor{ggbcolor26}{RGB}{68, 68, 68}
\definecolor{ggbcolor27}{RGB}{68, 68, 68}
\definecolor{ggbcolor28}{RGB}{153, 51, 0}
\definecolor{ggbcolor29}{RGB}{153, 51, 0}
\definecolor{ggbcolor30}{RGB}{153, 51, 0}
\definecolor{ggbcolor31}{RGB}{153, 51, 0}
\definecolor{ggbcolor32}{RGB}{153, 51, 0}
\definecolor{ggbcolor33}{RGB}{68, 68, 68}
\definecolor{ggbcolor34}{RGB}{68, 68, 68}
\definecolor{ggbcolor35}{RGB}{153, 51, 0}
\definecolor{ggbcolor36}{RGB}{153, 51, 0}
\definecolor{ggbcolor37}{RGB}{130, 0, 130}
\definecolor{ggbcolor38}{RGB}{130, 0, 130}
\definecolor{ggbcolor39}{RGB}{130, 0, 130}
\definecolor{ggbcolor40}{RGB}{130, 0, 130}
\definecolor{ggbcolor41}{RGB}{68, 68, 68}
\definecolor{ggbcolor42}{RGB}{68, 68, 68}
\definecolor{ggbcolor43}{RGB}{130, 0, 130}
\definecolor{ggbcolor44}{RGB}{130, 0, 130}
\definecolor{ggbcolor45}{RGB}{130, 0, 130}
\definecolor{ggbcolor46}{RGB}{68, 68, 68}
\definecolor{ggbcolor47}{RGB}{220, 0, 0}

\begin{longtable}[]{rlcp{4cm}p{3cm}}
\caption{Construction protocol for Problem 15}\label{GD15cp}\\
\toprule
No.&Name&Toolbar Icon&Description&Value\tabularnewline
\midrule
\endhead
\textcolor{ggbcolor2}{1}&\textcolor{ggbcolor2}{Polygon poly1}&\textcolor{ggbcolor2}{\raisebox{-1mm}{\includegraphics[width=0.5cm]{mode_regularpolygon}}}&\textcolor{ggbcolor2}{Polygon(A, B, 4)}&\textcolor{ggbcolor2}{$poly1 = 7.33$}\tabularnewline
\textcolor{ggbcolor3}{2}&\textcolor{ggbcolor3}{Segment f}&\textcolor{ggbcolor3}{\raisebox{-1mm}{\includegraphics[width=0.5cm]{mode_segment}}}&\textcolor{ggbcolor3}{Segment A, B}&\textcolor{ggbcolor3}{$f = 2.71$}\tabularnewline
\textcolor{ggbcolor9}{3}&\textcolor{ggbcolor9}{Polygon poly2}&\textcolor{ggbcolor9}{\raisebox{-1mm}{\includegraphics[width=0.5cm]{mode_regularpolygon}}}&\textcolor{ggbcolor9}{Polygon(D, C, 4)}&\textcolor{ggbcolor9}{$poly2 = 7.33$}\tabularnewline
\textcolor{ggbcolor16}{4}&\textcolor{ggbcolor16}{Polygon poly3}&\textcolor{ggbcolor16}{\raisebox{-1mm}{\includegraphics[width=0.5cm]{mode_regularpolygon}}}&\textcolor{ggbcolor16}{Polygon(C, B, 4)}&\textcolor{ggbcolor16}{$poly3 = 7.33$}\tabularnewline
\textcolor{ggbcolor23}{5}&\textcolor{ggbcolor23}{Polygon poly4}&\textcolor{ggbcolor23}{\raisebox{-1mm}{\includegraphics[width=0.5cm]{mode_regularpolygon}}}&\textcolor{ggbcolor23}{Polygon(A, D, 4)}&\textcolor{ggbcolor23}{$poly4 = 7.33$}\tabularnewline
\textcolor{ggbcolor30}{6}&\textcolor{ggbcolor30}{Polygon poly5}&\textcolor{ggbcolor30}{\raisebox{-1mm}{\includegraphics[width=0.5cm]{mode_regularpolygon}}}&\textcolor{ggbcolor30}{Polygon(B, A, 4)}&\textcolor{ggbcolor30}{$poly5 = 7.33$}\tabularnewline
\textcolor{ggbcolor37}{7}&\textcolor{ggbcolor37}{Segment g\_1}&\textcolor{ggbcolor37}{\raisebox{-1mm}{\includegraphics[width=0.5cm]{mode_segment}}}&\textcolor{ggbcolor37}{Segment E, K}&\textcolor{ggbcolor37}{$g_1 = 8.56$}\tabularnewline
\textcolor{ggbcolor38}{8}&\textcolor{ggbcolor38}{Polygon poly6}&\textcolor{ggbcolor38}{\raisebox{-1mm}{\includegraphics[width=0.5cm]{mode_regularpolygon}}}&\textcolor{ggbcolor38}{Polygon(I, D, 4)}&\textcolor{ggbcolor38}{$poly6 = 7.33$}\tabularnewline
\textcolor{ggbcolor45}{9}&\textcolor{ggbcolor45}{Segment l\_1}&\textcolor{ggbcolor45}{\raisebox{-1mm}{\includegraphics[width=0.5cm]{mode_segment}}}&\textcolor{ggbcolor45}{Segment N, H}&\textcolor{ggbcolor45}{$l_1 = 8.56$}\tabularnewline
\textcolor{ggbcolor46}{10}&\textcolor{ggbcolor46}{Point O}&\textcolor{ggbcolor46}{\raisebox{-1mm}{\includegraphics[width=0.5cm]{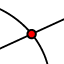}}}&\textcolor{ggbcolor46}{Intersection of g\_1 and l\_1}&\textcolor{ggbcolor46}{$O\, = \,\left(1.18, 3.29 \right)$}\tabularnewline
\textcolor{ggbcolor47}{11}&\textcolor{ggbcolor47}{Segment m\_1}&\textcolor{ggbcolor47}{\raisebox{-1mm}{\includegraphics[width=0.5cm]{mode_segment}}}&\textcolor{ggbcolor47}{Segment O, E}&\textcolor{ggbcolor47}{$m_1 = 1.71$}\tabularnewline
\bottomrule
\end{longtable}

\definecolor{ggbcolor0}{RGB}{40, 40, 170}
\definecolor{ggbcolor1}{RGB}{40, 40, 170}
\definecolor{ggbcolor2}{RGB}{40, 40, 170}
\definecolor{ggbcolor4}{RGB}{153, 51, 0}
\definecolor{ggbcolor5}{RGB}{153, 51, 0}
\definecolor{ggbcolor6}{RGB}{153, 51, 0}
\definecolor{ggbcolor7}{RGB}{68, 68, 68}
\definecolor{ggbcolor8}{RGB}{68, 68, 68}
\definecolor{ggbcolor9}{RGB}{153, 51, 0}
\definecolor{ggbcolor10}{RGB}{153, 51, 0}
\definecolor{ggbcolor11}{RGB}{130, 0, 130}
\definecolor{ggbcolor12}{RGB}{68, 68, 68}
\definecolor{ggbcolor14}{RGB}{68, 68, 68}
\definecolor{ggbcolor16}{RGB}{68, 68, 68}
\definecolor{ggbcolor17}{RGB}{220, 0, 0}

\begin{longtable}[]{rlcp{5cm}p{2.8cm}}
\caption{Construction protocol for Problem 23}\label{GD23cp}\\
\toprule
No.&Name&T.~Icon&Description&Value\tabularnewline
\midrule
\endhead
\textcolor{ggbcolor0}{1}&\textcolor{ggbcolor0}{Point A}&\textcolor{ggbcolor0}{\raisebox{-1mm}{\includegraphics[width=0.5cm]{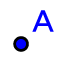}}}& &\textcolor{ggbcolor0}{$A\, = \,\left(0, 0 \right)$}\tabularnewline
\textcolor{ggbcolor1}{2}&\textcolor{ggbcolor1}{Point B}&\textcolor{ggbcolor1}{\raisebox{-1mm}{\includegraphics[width=0.5cm]{mode_point}}}& &\textcolor{ggbcolor1}{$B\, = \,\left(1, 0 \right)$}\tabularnewline
\textcolor{ggbcolor2}{3}&\textcolor{ggbcolor2}{Point B'}&\textcolor{ggbcolor2}{\raisebox{-1mm}{\includegraphics[width=0.5cm]{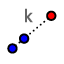}}}&\textcolor{ggbcolor2}{B dilated by factor 23 from A}&\textcolor{ggbcolor2}{$B'\, = \,\left(23, 0 \right)$}\tabularnewline
4&Segment n&\raisebox{-1mm}{\includegraphics[width=0.5cm]{mode_segment}}&Segment A, B'&$n = 23$\tabularnewline
\textcolor{ggbcolor4}{5}&\textcolor{ggbcolor4}{Polygon poly1}&\textcolor{ggbcolor4}{\raisebox{-1mm}{\includegraphics[width=0.5cm]{mode_regularpolygon}}}&\textcolor{ggbcolor4}{Polygon(B, A, 4)}&\textcolor{ggbcolor4}{$poly1 = 1$}\tabularnewline
\textcolor{ggbcolor5}{6}&\textcolor{ggbcolor5}{Segment f}&\textcolor{ggbcolor5}{\raisebox{-1mm}{\includegraphics[width=0.5cm]{mode_segment}}}&\textcolor{ggbcolor5}{Segment B, A}&\textcolor{ggbcolor5}{$f = 1$}\tabularnewline
\textcolor{ggbcolor11}{7}&\textcolor{ggbcolor11}{Line j}&\textcolor{ggbcolor11}{\raisebox{-1mm}{\includegraphics[width=0.5cm]{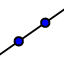}}}&\textcolor{ggbcolor11}{Line D, A}&\textcolor{ggbcolor11}{$j: -x - 1y\, = \,0$}\tabularnewline
\textcolor{ggbcolor12}{8}&\textcolor{ggbcolor12}{Point C'}&\textcolor{ggbcolor12}{\raisebox{-1mm}{\includegraphics[width=0.5cm]{mode_dilatefrompoint}}}&\textcolor{ggbcolor12}{C dilated by factor 11 from A}&\textcolor{ggbcolor12}{$C'\, = \,\left(0, -11 \right)$}\tabularnewline
9&Segment k&\raisebox{-1mm}{\includegraphics[width=0.5cm]{mode_segment}}&Segment C', B'&$k = 25.5$\tabularnewline
\textcolor{ggbcolor14}{10}&\textcolor{ggbcolor14}{Point E}&\textcolor{ggbcolor14}{\raisebox{-1mm}{\includegraphics[width=0.5cm]{mode_intersect}}}&\textcolor{ggbcolor14}{Intersection of j and k}&\textcolor{ggbcolor14}{$E\, = \,\left(7.44, -7.44 \right)$}\tabularnewline
11&Line l&\raisebox{-1mm}{\includegraphics[width=0.5cm]{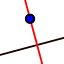}}&Line through E perpendicular to n&$l: x\, = \,7.44$\tabularnewline
\textcolor{ggbcolor16}{12}&\textcolor{ggbcolor16}{Point F}&\textcolor{ggbcolor16}{\raisebox{-1mm}{\includegraphics[width=0.5cm]{mode_intersect}}}&\textcolor{ggbcolor16}{Intersection of l and n}&\textcolor{ggbcolor16}{$F\, = \,\left(7.44, 0 \right)$}\tabularnewline
\textcolor{ggbcolor17}{13}&\textcolor{ggbcolor17}{Segment m}&\textcolor{ggbcolor17}{\raisebox{-1mm}{\includegraphics[width=0.5cm]{mode_segment}}}&\textcolor{ggbcolor17}{Segment F, E}&\textcolor{ggbcolor17}{$m = 7.44$}\tabularnewline
\bottomrule
\end{longtable}

\definecolor{ggbcolor0}{RGB}{27, 27, 205}
\definecolor{ggbcolor1}{RGB}{27, 27, 205}
\definecolor{ggbcolor3}{RGB}{27, 27, 205}
\definecolor{ggbcolor4}{RGB}{27, 27, 205}
\definecolor{ggbcolor6}{RGB}{68, 68, 68}
\definecolor{ggbcolor8}{RGB}{27, 27, 205}
\definecolor{ggbcolor11}{RGB}{220, 0, 0}
\definecolor{ggbcolor12}{RGB}{68, 68, 68}
\definecolor{ggbcolor13}{RGB}{0, 0, 255}

\begin{longtable}[]{rlcp{5cm}p{3cm}}
\caption{Construction protocol for Problem 47}\label{GD47cp}\\
\toprule
No.&Name&T.~Icon&Description&Value\tabularnewline
\midrule
\endhead
\textcolor{ggbcolor0}{1}&\textcolor{ggbcolor0}{Point B}&\textcolor{ggbcolor0}{\raisebox{-1mm}{\includegraphics[width=0.5cm]{mode_point}}}& &\textcolor{ggbcolor0}{$B\, = \,\left(-5.41, -5.8 \right)$}\tabularnewline
\textcolor{ggbcolor1}{2}&\textcolor{ggbcolor1}{Point A}&\textcolor{ggbcolor1}{\raisebox{-1mm}{\includegraphics[width=0.5cm]{mode_point}}}& &\textcolor{ggbcolor1}{$A\, = \,\left(-0.78, 4.18 \right)$}\tabularnewline
3&Line g&\raisebox{-1mm}{\includegraphics[width=0.5cm]{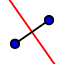}}&Perpendicular Bisector of AB&$g: 4.63x + 9.98y\, = \,-22.4$\tabularnewline
\textcolor{ggbcolor3}{4}&\textcolor{ggbcolor3}{Point D}&\textcolor{ggbcolor3}{\raisebox{-1mm}{\includegraphics[width=0.5cm]{mode_dilatefrompoint}}}&\textcolor{ggbcolor3}{B dilated by factor 8/11 from A}&\textcolor{ggbcolor3}{$D\, = \,\left(-4.15, -3.08 \right)$}\tabularnewline
\textcolor{ggbcolor4}{5}&\textcolor{ggbcolor4}{Point A'}&\textcolor{ggbcolor4}{\raisebox{-1mm}{\includegraphics[width=0.5cm]{mode_dilatefrompoint}}}&\textcolor{ggbcolor4}{A dilated by factor 1/8 from D}&\textcolor{ggbcolor4}{$A'\, = \,\left(-3.73, -2.17 \right)$}\tabularnewline
6&Circle c&\raisebox{-1mm}{\includegraphics[width=0.5cm]{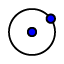}}&Circle through A' with center A&$c: (x + 0.78)^{2} + (y - 4.18)^{2} = 49.04$\tabularnewline
\textcolor{ggbcolor6}{7}&\textcolor{ggbcolor6}{Point O}&\textcolor{ggbcolor6}{\raisebox{-1mm}{\includegraphics[width=0.5cm]{mode_intersect}}}&\textcolor{ggbcolor6}{Intersection of c and g}&\textcolor{ggbcolor6}{$O\, = \,\left(0.84, -2.63 \right)$}\tabularnewline
8&Circle d&\raisebox{-1mm}{\includegraphics[width=0.5cm]{mode_circle2}}&Circle through A with center O&$d: (x - 0.84)^{2} + (y + 2.63)^{2} = 49.04$\tabularnewline
\textcolor{ggbcolor8}{9}&\textcolor{ggbcolor8}{Point E}&\textcolor{ggbcolor8}{\raisebox{-1mm}{\includegraphics[width=0.5cm]{mode_mirroratpoint}}}&\textcolor{ggbcolor8}{D mirrored at O}&\textcolor{ggbcolor8}{$E\, = \,\left(5.82, -2.19 \right)$}\tabularnewline
10&Line h&\raisebox{-1mm}{\includegraphics[width=0.5cm]{mode_join}}&Line A, E&$h: 6.37x + 6.6y\, = \,22.63$\tabularnewline
11&Segment f&\raisebox{-1mm}{\includegraphics[width=0.5cm]{mode_segment}}&Segment A, B&$f = 11$\tabularnewline
\textcolor{ggbcolor11}{12}&\textcolor{ggbcolor11}{Segment i}&\textcolor{ggbcolor11}{\raisebox{-1mm}{\includegraphics[width=0.5cm]{mode_segment}}}&\textcolor{ggbcolor11}{Segment D, A'}&\textcolor{ggbcolor11}{$i = 1$}\tabularnewline
\textcolor{ggbcolor12}{13}&\textcolor{ggbcolor12}{Point C}&\textcolor{ggbcolor12}{\raisebox{-1mm}{\includegraphics[width=0.5cm]{mode_intersect}}}&\textcolor{ggbcolor12}{Intersection of d and h}&\textcolor{ggbcolor12}{$C\, = \,\left(7.7, -4.01 \right)$}\tabularnewline
\textcolor{ggbcolor13}{14}&\textcolor{ggbcolor13}{Segment j}&\textcolor{ggbcolor13}{\raisebox{-1mm}{\includegraphics[width=0.5cm]{mode_segment}}}&\textcolor{ggbcolor13}{Segment C, E}&\textcolor{ggbcolor13}{$j = 2.62$}\tabularnewline
\bottomrule
\end{longtable}

\end{document}